\documentclass[11pt,reqno]{amsart}

\usepackage[T1]{fontenc}
\usepackage{lmodern}
\usepackage{amsmath,amssymb,amsthm,mathtools}
\usepackage[margin=1.12in]{geometry}
\usepackage{microtype}
\usepackage{enumitem}
\usepackage{booktabs}
\usepackage{graphicx}
\usepackage[dvipsnames]{xcolor}
\usepackage[colorlinks=true,linkcolor=MidnightBlue,citecolor=MidnightBlue,urlcolor=MidnightBlue]{hyperref}
\usepackage[capitalise,nameinlink]{cleveref}
\hypersetup{pdftitle={Forbidden subset sums in a linear range: an inverse theorem and the Erd\H{o}s--Graham irregularity problem},pdfauthor={Mateo Matijasevick}}

\numberwithin{equation}{section}
\newtheorem{theorem}{Theorem}[section]
\newtheorem{proposition}[theorem]{Proposition}
\newtheorem{lemma}[theorem]{Lemma}
\newtheorem{corollary}[theorem]{Corollary}
\newtheorem{conjecture}[theorem]{Conjecture}
\theoremstyle{definition}

\newtheorem{problem}[theorem]{Problem}
\theoremstyle{remark}
\newtheorem{remark}[theorem]{Remark}
\crefname{problem}{Problem}{Problems}
\crefname{conjecture}{Conjecture}{Conjectures}

\newcommand{\N}{\mathbb{N}}
\newcommand{\Z}{\mathbb{Z}}
\newcommand{\eps}{\varepsilon}
\newcommand{\card}[1]{\lvert #1\rvert}
\newcommand{\floor}[1]{\lfloor #1\rfloor}
\newcommand{\ceil}[1]{\lceil #1\rceil}
\newcommand{\snd}{\operatorname{snd}}
\newcommand{\lcmop}{\operatorname{lcm}}
\newcommand{\Sig}{\Sigma}
\newcommand{\lean}[1]{\texttt{\detokenize{#1}}}

\title[Forbidden subset sums in a linear range]{Forbidden subset sums in a linear range:\\ an inverse theorem and the Erd\H{o}s--Graham irregularity problem}

\author{Mateo Matijasevick}
\address{Purdue University, West Lafayette, Indiana, USA}
\email{mmatijas@purdue.edu}

\subjclass[2020]{Primary 11B13; Secondary 11B30, 11P70, 68V20}
\keywords{Subset sums, Erd\H{o}s problems, inverse additive problems, least non-divisor, Freiman's $3k-4$ theorem, zero-sum theorems, formal verification, Lean}

\date{}

\begin{document}

\begin{abstract}
For $c>0$ let $f_c(n)$ be the largest size of a set $A\subseteq\{1,\dots,\lfloor cn\rfloor\}$ such that $n$ is not a sum of distinct elements of $A$. Erd\H{o}s and Graham asked whether $f_c(n)$ depends on $n$ in an irregular way (Erd\H{o}s Problem \#361). We prove that $f_c(n)=\lfloor cn\rfloor-\lceil n/2\rceil$ for $c\ge 1$, while for $0<c<1$ the ratio $f_c(n)/n$ has $\limsup$ equal to $c/2$, attained along odd $n$, and, along the integers $n$ whose least non-divisor equals a prime power $q$, tends to
\[
\max\Bigl\{\frac cq,\; c-\frac12+\frac{1-c}{q}\Bigr\}
\qquad\text{whenever } 2^{q-2}\min\{c,1-c\}\le 1 .
\]
Consequently $f_c(n)/n$ converges if and only if $c\ge1$; for $0<c<1$ it has at least two limit points, and arbitrarily many as $c\to0^+$ or $c\to1^-$. For $c=1/2$ the set of limit points is exactly $\{1/6,1/4\}$. The upper bounds follow from an inverse theorem of independent interest: if $s\ge1$, $q$ is the least positive integer not dividing $s$, and $\eps>0$, then every $A\subseteq\{1,\dots,E\}$ with $|A|\ge E/q+\eps T$ has a subset summing to $sT$, uniformly for all sufficiently large $T$ and all $E\le sT/2^{q-2}$. For $s=2$ this extends Alon's theorem on subsets of $\{1,\dots,N\}$ with no subset sum equal to $2N$. The results for $c\ge1$ and the $\limsup$, as well as the limits and the inverse theorem for $q\in\{2,3\}$, have been formalized and machine-checked in the Lean~4 proof assistant.
\end{abstract}

\maketitle

\section{Introduction}\label{sec:intro}

\subsection{The problem}
Throughout, $\N=\{1,2,3,\dots\}$, and for integers $a\le b$ we write $[a,b]=\{a,a+1,\dots,b\}$; if $a>b$ then $[a,b]=\emptyset$. For a finite set $A\subseteq\N$ let
\[
 \Sig(A)=\Bigl\{\sum_{b\in B}b:\ \emptyset\neq B\subseteq A\Bigr\}
\]
be the set of its subset sums; we say that $A$ \emph{avoids} $n$ if $n\notin\Sig(A)$. For integers $M\ge0$ and $n\ge1$ put
\[
 F(M,n)=\max\bigl\{\card A:\ A\subseteq[1,M],\ n\notin\Sig(A)\bigr\},
 \qquad
 f_c(n)=F(\floor{cn},n)\quad(c>0).
\]
Erd\H{o}s and Graham~\cite{EG80} asked the following question, which is Problem~\#361 in Bloom's database~\cite{Bl} (quoted from there; it is listed as open at the time of writing): \emph{Let $c>0$ and $n$ be some large integer. What is the size of the largest $A\subseteq\{1,\dots,\floor{cn}\}$ such that $n$ is not a sum of a subset of $A$? Does this depend on $n$ in an irregular way?}

Two elementary observations frame the problem. For a positive integer $m$ let $\snd(m)$ denote the least positive integer that does not divide $m$. If $q\ge2$ does not divide $n$, then every subset sum of the multiples of $q$ is a multiple of $q$; hence
\begin{equation}\label{eq:trivial-lower}
 f_c(n)\ \ge\ \Bigl\lfloor\frac{\floor{cn}}{\snd(n)}\Bigr\rfloor .
\end{equation}
Secondly, intervals avoid $n$: for every $k\ge1$ the integers in $(n/(k+1),n/k)$ avoid $n$, since a sum of at most $k$ of them is smaller than $n$ and a sum of at least $k+1$ of them is larger. The first observation suggests that $f_c(n)$ should depend on the arithmetic of $n$, but a proof of irregularity requires \emph{upper} bounds that are sensitive to that arithmetic. Providing such bounds is the purpose of this paper.

\subsection{Results}
The case $c\ge1$ is elementary; we record it for completeness (it has also been noted in the discussion of the problem on~\cite{Bl}).

\begin{theorem}[Regular range]\label{thm:cge1}
For all integers $1\le n\le M$ we have $F(M,n)=M-\ceil{n/2}$. Consequently, for every $c\ge1$ and $n\ge1$,
\[
 f_c(n)=\floor{cn}-\ceil{n/2},\qquad\text{and}\qquad \frac{f_c(n)}{n}\longrightarrow c-\frac12 .
\]
\end{theorem}

The range $0<c<1$ is governed by the following inverse theorem.

\begin{theorem}[Inverse theorem]\label{thm:inverse}
Let $s\ge1$ be an integer, let $q=\snd(s)$ and let $\eps>0$. There is $T_0=T_0(s,\eps)$ with the following property. Let $T\ge T_0$ and $E\ge1$ be integers with
\begin{equation}\label{eq:range}
 2^{q-2}E\le sT .
\end{equation}
If $A\subseteq[1,E]$ satisfies $\card A\ge E/q+\eps T$, then $sT$ is a sum of at most $2^{q-2}(\ceil{2s/\eps}+1)$ distinct elements of $A$; in particular $sT\in\Sig(A)$.
\end{theorem}

The density $1/q$ is optimal whenever $q\nmid sT$, by the multiples of $q$ in $[1,E]$; and for $q\in\{2,3\}$ the range~\eqref{eq:range} is optimal as well (\cref{prop:sharp}). For $s=2$ and $E=T$, \cref{thm:inverse} contains the upper bound in Alon's theorem~\cite[Corollary~2.6]{Al87} that a subset of $[1,N]$ with no subset sum equal to $2N$ has at most $(1/3+o(1))N$ elements. With $E=\floor{xT}$, \cref{thm:inverse} takes the following form.

\begin{corollary}\label{cor:interval}
Let $s\ge1$, $q=\snd(s)$, $\eps>0$, and $0<x\le s/2^{q-2}$ (for instance, any $0<x\le1$). For all sufficiently large $T$, every $A\subseteq[1,\floor{xT}]$ with $\card A>(x/q+\eps)T$ satisfies $sT\in\Sig(A)$. Moreover, if $q=p^a$ with $p$ prime, then
\[
 F(\floor{xT},sT)=\Bigl(\frac xq+o(1)\Bigr)T\qquad\text{as } T\to\infty,\ p\nmid T .
\]
\end{corollary}

For $0<c<1$ let
\[
 \varphi_q(c)=\max\Bigl\{\frac cq,\ c-\frac12+\frac{1-c}q\Bigr\}
 =\begin{cases} c/q, & c\le 1/2,\\ c-\tfrac12+\tfrac{1-c}q, & c\ge1/2,\end{cases}
\]
and let $\mathcal Q(c)$ be the set of prime powers $q$ with $2^{q-2}\min\{c,1-c\}\le1$. Note that $\{2,3\}\subseteq\mathcal Q(c)$ for every $c\in(0,1)$, that $4\in\mathcal Q(c)$ if and only if $\min\{c,1-c\}\le1/4$, and so on (see \cref{tab:classes}). Recall that $\snd(n)$ is always a prime power (\cref{lem:snd}).

\begin{theorem}[Asymptotics for $0<c<1$]\label{thm:main}
Let $0<c\le1$.
\begin{enumerate}[label=\textup{(\roman*)}]
\item $f_c(n)\le(c/2+o(1))\,n$ as $n\to\infty$, and $f_c(n)\ge cn/2-1$ for every odd $n$. Hence
\[
 \limsup_{n\to\infty}\frac{f_c(n)}n=\frac c2=\lim_{\substack{n\to\infty\\ n\ \mathrm{odd}}}\frac{f_c(n)}n .
\]
\item Assume $c<1$ and let $q\in\mathcal Q(c)$. Then $f_c(n)\le(\varphi_q(c)+o(1))\,n$ as $n\to\infty$ through the integers with $\snd(n)\ge q$, and $f_c(n)\ge\varphi_q(c)\,n-q$ whenever $\snd(n)=q$. In particular
\[
 \lim_{\substack{n\to\infty\\ \snd(n)=q}}\frac{f_c(n)}n=\varphi_q(c)
 \qquad\text{and}\qquad
 \limsup_{\substack{n\to\infty\\ \snd(n)\ge q}}\frac{f_c(n)}n=\varphi_q(c).
\]
\end{enumerate}
\end{theorem}

\begin{table}[ht]
\centering
\caption{The first classes of \cref{thm:main}(ii). The class $\{n:\snd(n)=q\}$ consists of the multiples of $\lcmop(1,\dots,q-1)$ that are not multiples of $q$.}\label{tab:classes}
\begin{tabular}{@{}cllll@{}}
\toprule
$q$ & $\{n:\snd(n)=q\}$ & condition $q\in\mathcal Q(c)$ & $\varphi_q(c)$, $c\le\frac12$ & $\varphi_q(c)$, $c\ge\frac12$\\
\midrule
$2$ & $n$ odd & always & $c/2$ & $c/2$\\
$3$ & $n\equiv 2,4 \pmod 6$ & always & $c/3$ & $2c/3-1/6$\\
$4$ & $n\equiv 6\pmod{12}$ & $c\le\frac14$ or $c\ge\frac34$ & $c/4$ & $3c/4-1/4$\\
$5$ & $n\equiv 12,24,36,48\pmod{60}$ & $c\le\frac18$ or $c\ge\frac78$ & $c/5$ & $4c/5-3/10$\\
$7$ & $n\equiv 60k\pmod{420}$, $7\nmid k$ & $c\le\frac1{32}$ or $c\ge\frac{31}{32}$ & $c/7$ & $6c/7-5/14$\\
\bottomrule
\end{tabular}
\end{table}

\begin{corollary}[Irregularity]\label{cor:irregular}
Let $0<c<1$. The numbers $\varphi_q(c)$, $q\in\mathcal Q(c)$, are pairwise distinct limit points of $f_c(n)/n$. In particular $f_c(n)/n$ does not converge, and the number of its limit points is at least the number of prime powers not exceeding $2+\log_2(1/\min\{c,1-c\})$, which tends to infinity as $c\to0^+$ or $c\to1^-$. Together with \cref{thm:cge1}: $f_c(n)/n$ converges if and only if $c\ge1$.
\end{corollary}

For one value of $c$ the picture is complete.

\begin{corollary}[$c=1/2$]\label{cor:half}
$f_{1/2}(n)/n\to1/4$ along odd $n$ and $f_{1/2}(n)/n\to1/6$ along even $n$. Hence the set of limit points of $f_{1/2}(n)/n$ is exactly $\{1/6,1/4\}$.
\end{corollary}

The hypothesis $q\in\mathcal Q(c)$ in \cref{thm:main}(ii) cannot simply be dropped: for $2/3\le c<3/4$ an interval construction gives $f_c(n)\ge n/4-1$ for every $n$, which exceeds $\varphi_q(c)\,n$ for all large $q$ (\cref{rem:intervals}). Determining the behaviour of $f_c(n)/n$ along the classes $\snd(n)=q\notin\mathcal Q(c)$ remains open (\cref{sec:open}).

\subsection{Relation to previous work}\label{subsec:literature}
In Alon's notation~\cite{Al87}, $f(N,m)$ is the largest size of $A\subseteq[1,N]$ with $m\notin\Sig(A)$; thus $f_c(n)=f(\floor{cn},n)$, and Problem~\#361 concerns the \emph{linear regime} $m\asymp N$. Alon introduced the function $\snd$, observed the lower bound $f(N,m)\ge\floor{N/\snd(m)}$ (the multiples of $\snd(m)$), and proved~\cite[Theorem~1.2]{Al87} that $f(N,m)\le C(\eps)\,N/\snd(m)$ for $N^{1+\eps}\le m\le N^2/\log^2N$. He conjectured~\cite[Conjecture~4.3]{Al87} that $f(N,m)=(1+o(1))N/\snd(m)$ in a polynomial range of $m$ (see~\cite{TVW09} for the precise formulation). This was proved by Lipkin~\cite{Li89} for $cN\log^6N<m<N^{3/2}/\log^3N$; Alon and Freiman~\cite{AF88} determined $f(N,m)=\floor{N/\snd(m)}+\snd(m)-2$ exactly for $3N^{5/3+\eps}<m<N^2/(20\log^2N)$; Tran, Vu and Wood~\cite{TVW09} established the asymptotic formula for $cN(\log N)^{1+\eps}\le m\le N^2/(9\log^2N)$; and Conlon, Fox and Pham~\cite{CFP} determined $f(N,m)$ exactly for $CN\log N\le m\le N^2/(12\log^2N)$ and asymptotically for $N^2/(12\log^2N)\le m\le\binom{N+1}2$, and showed that $f(N,m)$ exceeds $N/\snd(m)$ by a factor close to $2$ for some $m$ with $8N<m<N(\log N)/8$. None of these results applies when $m\asymp N$. In this regime the only previous asymptotic result we are aware of is Alon's $f(N,2N)=(1/3+o(1))N$~\cite[Corollary~2.6]{Al87}, which corresponds to $c=1/2$ and even $n$. \Cref{thm:main} shows that for $c\le1/2$ the formula $f(N,m)\sim N/\snd(m)$ (with $N=\floor{cn}$, $m=n$) holds along the classes $\snd(n)=q\in\mathcal Q(c)$, while for $c>1/2$ and $q\ge3$ it fails: the extremal sets are then not the multiples of $\snd(n)$ but the mixed sets of \cref{prop:constructions}(b).

Our arguments combine three classical ingredients: Roth's theorem~\cite{Ro53}, used, as in Alon's proof of~\cite[Proposition~2.5]{Al87}, to pass from sums with repetitions to sums of distinct elements (an analogue for sets of positive density, via Szemer\'edi's theorem, is~\cite[Lemmas~1 and~2]{ENS88}); an elementary ``central interval'' lemma for $2B$ and $hB$; and Freiman's $(3k-4)$-theorem~\cite{Fr73} (see also~\cite{LS95,HP02,Na96}), iterated along the dyadic sumsets $2^jB$. Related results on long intervals and arithmetic progressions in subset sums of dense sets include~\cite{EF90,NS89,Le97,Le03,SV06}. The zero-sum route of \cref{sec:zerosum} uses Alon's bounded zero-sum theorem~\cite[Theorem~1.1]{Al87} or, for prime moduli, the Dias da Silva--Hamidoune theorem~\cite{DSH94,ANR96}.

The constructions in \cref{prop:constructions} are simple, and several of them were proposed in the discussion of Problem~\#361 on~\cite{Bl}: the multiples of the least non-divisor, the mixed set of multiples of $3$ together with a residue class above $n/2$ (for $2\mid n$, $3\nmid n$), intervals of the form $(n/(k+1),n/k)$, and the formula for $c\ge1$; the numerical values reported there for $c=3/4$ and $100\le n\le104$ are reproduced in \cref{sec:numerics}. The contribution of the present paper is the matching upper bounds, i.e.\ \cref{thm:inverse} and its consequences.

\subsection{Formal verification}
The following statements have been formalized in Lean~4 with Mathlib~\cite{Lean4,mathlib} and machine-checked: \cref{thm:cge1}, both the formula for $F(M,n)$ and the limit $f_c(n)/n\to c-1/2$; the representation statement of \cref{thm:inverse} for $q\in\{2,3\}$, that is, for all $s$ not divisible by $6$; \cref{thm:main}(i), as the upper bound $f_c(n)\le(c/2+\eps)n$ for large $n$ together with the limit $c/2$ along odd $n$; the limits of \cref{thm:main}(ii) along $\snd(n)=q$ for $q\in\{2,3\}$; and the non-convergence of $f_c(n)/n$ for every $0<c<1$, once through the two limit points given by $q=2$ and $q=3$ and once through the zero-sum argument of \cref{sec:zerosum}. Each formal statement depends only on the three standard axioms of Lean's logic (propositional extensionality, soundness of quotients, and choice). See \cref{sec:lean}.

\subsection{Organization}
\Cref{sec:prelim} collects the lemmas. \Cref{sec:inverse} proves \cref{thm:inverse,cor:interval} and discusses sharpness. \Cref{sec:cge1} proves \cref{thm:cge1}, and \cref{sec:asymp} proves \cref{thm:main} and \cref{cor:irregular,cor:half}. \Cref{sec:zerosum} gives the zero-sum proof of non-convergence. \Cref{sec:numerics} reports exact computations, \cref{sec:lean} describes the formalization, and \cref{sec:open} lists open problems.

\section{Preliminaries}\label{sec:prelim}

\subsection{Notation}
For $B\subseteq\Z$ and an integer $h\ge1$ let $hB=\{b_1+\dots+b_h:\ b_i\in B\}$ (repetitions allowed), and $0B=\{0\}$. If $0\in B$ then $hB\subseteq(h+1)B$, and $h(kB)=(hk)B$. For a nonempty finite $G\subseteq\N$, $\gcd(G)$ is the greatest common divisor of its elements, and $\sigma(X)$ denotes the sum of the elements of a finite set or multiset $X$. Let $r_3(N)$ be the largest size of a subset of $[1,N]$ containing no three-term arithmetic progression $x<y<z$. By Roth's theorem~\cite{Ro53}, $r_3(N)=o(N)$; quantitative bounds (e.g.~\cite{KM23}) make every $o(\cdot)$ below effective. Finally, $\Lambda_k=\lcmop(1,\dots,k)$ for $k\ge1$.

\begin{lemma}[Least non-divisors]\label{lem:snd}
Let $m\ge1$ and $q=\snd(m)$.
\begin{enumerate}[label=\textup{(\alph*)}]
\item $q\ge2$ and $q$ is a prime power.
\item $\Lambda_{q-1}\mid m$, and $\snd(\Lambda_{q-1})=q$.
\item If $q=p^a$ with $p$ prime, then $v_p(m)=a-1$; consequently, for every $T\in\N$, $q\mid mT$ if and only if $p\mid T$.
\item $\Lambda_k\ge2^{k-1}$ for every $k\ge1$; consequently $m\ge\Lambda_{q-1}\ge2^{q-2}$.
\end{enumerate}
\end{lemma}

\begin{proof}
(a) Since $1\mid m$, $q\ge2$. If $q$ were not a prime power, we could write $q=ab$ with $\gcd(a,b)=1$ and $2\le a,b<q$; then $a\mid m$ and $b\mid m$, hence $ab=q\mid m$, a contradiction.

(b) Every $k\in[1,q-1]$ divides $m$, so $\Lambda_{q-1}\mid m$. Every $k\le q-1$ divides $\Lambda_{q-1}$, while $q=p^a\nmid\Lambda_{q-1}$ because $v_p(\Lambda_{q-1})=\max\{j:p^j\le q-1\}=a-1$.

(c) If $a\ge2$ then $p^{a-1}\le q-1$ divides $m$; hence $v_p(m)\ge a-1$ in all cases, and $v_p(m)\le a-1$ because $q\nmid m$. Thus $v_p(mT)=a-1+v_p(T)$, which is at least $a$ if and only if $p\mid T$.

(d) For $1\le j\le k$,
\[
 \frac1{j\binom kj}=\int_0^1x^{j-1}(1-x)^{k-j}\,dx=\sum_{i=0}^{k-j}(-1)^i\binom{k-j}{i}\frac1{j+i},
\]
and every denominator $j+i$ lies in $[1,k]$; hence $j\binom kj$ divides $\Lambda_k$ (this is Nair's argument~\cite{Na82}). Since $j\binom kj=k\binom{k-1}{j-1}$ and $\max_j\binom{k-1}{j-1}\ge2^{k-1}/k$, we obtain $\Lambda_k\ge2^{k-1}$. Finally $m\ge\Lambda_{q-1}$ by (b), and $\Lambda_{q-1}\ge2^{q-2}$.
\end{proof}

\subsection{Good cores}
For a finite $A\subseteq\N$ and $g\in A$ put
\[
 \nu_A(g)=\#\{u\in A:\ u<g,\ 2g-u\in A\},
\]
the number of nontrivial three-term progressions in $A$ with middle term $g$. For $m\ge1$ let $C_m(A)=\{g\in A:\nu_A(g)<m\}$. The following lemma is a mild variant of the argument in Alon's proof of~\cite[Proposition~2.5]{Al87}.

\begin{lemma}[Distinctification]\label{lem:distinct}
Let $m\ge1$ and let $A\subseteq\N$ be finite.
\begin{enumerate}[label=\textup{(\alph*)}]
\item Let $G\subseteq A$ satisfy $\nu_A(g)\ge m$ for all $g\in G$. If $t=g_1+\dots+g_h$ with $1\le h\le m$ and $g_1,\dots,g_h\in G$ \textup{(}not necessarily distinct\textup{)}, then $t=a_1+\dots+a_h$ for some pairwise distinct $a_1,\dots,a_h\in A$. In particular $t\in\Sig(A)$.
\item If $A\subseteq[1,E]$, then $\card{C_m(A)}\le2\sqrt{2m}\;r_3(E)$.
\end{enumerate}
\end{lemma}

\begin{proof}
(a) We work with finite multisets $X$ of elements of $A$; $\card X$ counts multiplicity, $S(X)$ is the support and $e(X)=\card X-\card{S(X)}$ is the excess. Call $X$ \emph{admissible} if every element of multiplicity at least $2$ in $X$ belongs to $G$. The multiset $X_0=\{\!\{g_1,\dots,g_h\}\!\}$ is admissible, has cardinality $h$ and sum $t$.

Let $X$ be admissible with $\card X=h$, $\sigma(X)=t$ and $e(X)>0$, and let $g$ have multiplicity at least $2$ in $X$; then $g\in G$. Let $X'$ be obtained by deleting two copies of $g$; then $\card{S(X')}\le h-2\le m-2$. For each of the $\nu_A(g)\ge m$ elements $u\in A$ with $u<g$ and $2g-u\in A$ consider $P_u=\{u,2g-u\}$. These pairs are pairwise disjoint: their smaller elements are $<g$ and their larger elements are $>g$, and $u\mapsto u$, $u\mapsto2g-u$ are injective. Since every element of $S(X')$ lies in at most one pair, some $P_u$ is disjoint from $S(X')$. Let $X^+=X'+\{\!\{u,2g-u\}\!\}$. Then $\card{X^+}=h$ and $\sigma(X^+)=t$; the elements $u$ and $2g-u$ are distinct and have multiplicity $1$ in $X^+$, and no other multiplicity increased, so $X^+$ is admissible. Deleting two copies of $g$ lowers $\card X$ by $2$ and $\card{S(X)}$ by at most $1$, and adding $u,2g-u$ raises both by $2$; hence $e(X^+)\le e(X)-1$. Iterating, we reach an admissible multiset of cardinality $h$, sum $t$ and excess $0$, i.e.\ $h$ distinct elements of $A$ with sum $t$.

(b) Let $C=C_m(A)$ and $N=\card C$. If $m=1$, then $C$ contains no progression $u<g<2g-u$ (its middle term would have $\nu_A(g)\ge1$), so $N\le r_3(E)$. Let $m\ge2$, and let $\mathcal H$ be the $3$-uniform hypergraph on $C$ whose edges are the sets $\{u,g,2g-u\}\subseteq C$ with $u<g$. Each edge has a unique middle element $g$, and at most $\nu_A(g)\le m-1$ edges have middle element $g$; hence $\mathcal H$ has at most $(m-1)N$ edges. Choose $R\subseteq C$ at random, containing each vertex independently with probability $p=(2(m-1))^{-1/2}\in(0,1)$. The expected number of vertices of $R$ is $pN$, and the expected number of edges inside $R$ is at most $p^3(m-1)N=pN/2$. Fix an outcome with $\card R-\#\{\text{edges inside }R\}\ge pN/2$ and delete one vertex from each edge inside $R$. The resulting set $R'\subseteq[1,E]$ has at least $pN/2$ elements and contains no edge of $\mathcal H$, i.e.\ no three-term progression. Therefore $pN/2\le r_3(E)$, that is, $N\le2\sqrt{2(m-1)}\,r_3(E)\le2\sqrt{2m}\,r_3(E)$.
\end{proof}

\subsection{Sumsets of sets containing \texorpdfstring{$0$}{0}}

\begin{lemma}[Central interval]\label{lem:central}
Let $L\ge1$ and $B\subseteq[0,L]$ with $0,L\in B$. Put $N=\card B$ and suppose $L\le2N-4$. Let $\gamma=2L+2-2N$. Then $0\le\gamma\le L-2$, and for every integer $h\ge2$,
\[
 [\gamma,\ hL-\gamma]\subseteq hB .
\]
\end{lemma}

\begin{proof}
Since $N\le L+1$ we have $\gamma\ge0$, and $\gamma\le L-2$ is equivalent to $L+4\le2N$.

\emph{Step 1: $[\gamma,2L-\gamma]\subseteq2B$.} Let $\gamma\le z\le2L-\gamma$ and $U_z=\{u\in[0,L]:z-u\in[0,L]\}$. Then $\card{U_z}=\min\{z,2L-z,L\}+1\ge\gamma+1$. The set $[0,L]\setminus B$ has $L+1-N=\gamma/2$ elements, so at most $\gamma/2$ elements $u\in U_z$ satisfy $u\notin B$, and at most $\gamma/2$ satisfy $z-u\notin B$. Hence some $u\in U_z$ has $u\in B$ and $z-u\in B$.

\emph{Step 2.} Let $D=B\setminus\{L\}\subseteq[0,L-1]$; then $\card D=N-1$ and $2\card D=2N-2>L$. For $r\in[0,L-1]$, the sets $D$ and $r-D$, reduced modulo $L$, have $2N-2>L$ elements in total (reduction is injective on $[0,L-1]$), so they meet: there are $b_1,b_2\in D$ with $b_1+b_2\equiv r\pmod L$. Since $0\le b_1+b_2\le2L-2$, we get $b_1+b_2\in\{r,r+L\}$.

\emph{Step 3.} Let $h\ge2$ and $\gamma\le z\le hL-\gamma$. If $z\le2L-\gamma$, write $z=b+b'$ by Step~1 and add $h-2$ zeros. If $z\ge(h-2)L+\gamma$, then $z-(h-2)L\in[\gamma,2L-\gamma]$, and $z=b+b'+(h-2)L$ with $b,b'\in B$ by Step~1. Otherwise $2L-\gamma<z<(h-2)L+\gamma$; write $z=kL+r$ with $0\le r<L$. As $z>2L-\gamma\ge L+2$, we have $k\ge1$; as $z<(h-2)L+\gamma\le(h-1)L-2$, we have $k\le h-2$. Take $b_1,b_2$ from Step~2. If $b_1+b_2=r$ then $z=b_1+b_2+k\cdot L$ uses $k+2\le h$ summands; if $b_1+b_2=r+L$ then $z=b_1+b_2+(k-1)\cdot L$ uses $k+1\le h$ summands. Adding zeros gives exactly $h$ summands from $B$.
\end{proof}

We use Freiman's $(3k-4)$-theorem in the following normalized form; see~\cite{Fr73}, \cite[Chapter~1]{Na96} and~\cite{LS95}. (It is proved from scratch in the Lean development of \cref{sec:lean}.)

\begin{theorem}[Freiman]\label{thm:freiman}
Let $B\subseteq[0,L]$ with $0,L\in B$, $\gcd(B)=1$ and $k=\card B\ge3$. Then $\card{2B}\ge\min\{3k-3,\ k+L\}$. In particular, if $L\ge2k-3$ then $\card{2B}\ge3k-3$.
\end{theorem}

\begin{lemma}[Bounded additive basis]\label{lem:basis}
Let $q\ge2$ be an integer, let $\eta,K>0$, and put $H=2^{q-2}(\ceil K+1)$. Let $L$ be an integer with $L>3/(2\eta)$ and $L\ge3q$, and let $B\subseteq[0,L]$ satisfy $0,L\in B$, $\gcd(B)=1$ and $\card B\ge(1/q+\eta)L$. Then every integer $y$ with $2^{q-2}L\le y\le KL$ belongs to $HB$.
\end{lemma}

\begin{proof}
For $j\ge0$ let $B_j=2^jB$ and $N_j=\card{B_j}$. Since $0\in B$, we have $B\subseteq B_j\subseteq[0,2^jL]$ and $0,2^jL\in B_j$; in particular $\gcd(B_j)=1$ and $N_j\ge N_0\ge L/q\ge3$.

We claim that $2^jL\le2N_j-4$ for some $j\in[0,q-2]$. Suppose not, so that $2^jL\ge2N_j-3$ for all $0\le j\le q-2$. For $0\le j\le q-3$, \cref{thm:freiman} applied to $B_j$ gives $N_{j+1}=\card{2B_j}\ge3N_j-3$, i.e.\ $N_{j+1}-\frac32\ge3(N_j-\frac32)$. Hence $N_{q-2}-\frac32\ge3^{q-2}(N_0-\frac32)$, and
\[
 2^{q-2}L\ \ge\ 2N_{q-2}-3\ \ge\ 2\cdot3^{q-2}N_0-3^{q-1}\ \ge\ 2\cdot3^{q-2}\Bigl(\frac1q+\eta\Bigr)L-3^{q-1}.
\]
Now $(3/2)^{q-2}\ge q/2$ for all $q\ge2$ (equality for $q=2,3$; the left side grows by the factor $3/2\ge(q+1)/q$ when $q$ increases by one), so $2\cdot3^{q-2}/q\ge2^{q-2}$ and therefore $0\ge2\cdot3^{q-2}\eta L-3^{q-1}$, i.e.\ $L\le3/(2\eta)$, a contradiction. (For $q=2$ no application of \cref{thm:freiman} is needed.)

Fix such a $j$ and put $L_j=2^jL$. By \cref{lem:central} applied to $B_j$ (with $L_j\le2N_j-4$), the number $\gamma_j=2L_j+2-2N_j$ satisfies $\gamma_j\le L_j-2$, and $[\gamma_j,hL_j-\gamma_j]\subseteq hB_j$ for all $h\ge2$. Let $h=\ceil K+1\ge2$ and $2^{q-2}L\le y\le KL$. Then $y\ge2^{q-2}L\ge L_j\ge\gamma_j$ and $y\le KL\le(h-1)L_j\le hL_j-\gamma_j$. Hence $y\in hB_j=(h2^j)B\subseteq(h2^{q-2})B=HB$, using $0\in B$.
\end{proof}

\section{The inverse theorem}\label{sec:inverse}

\begin{proof}[Proof of \cref{thm:inverse}]
Let $s,\eps$ be given and $q=\snd(s)$. Put
\[
 K=\frac{2s}{\eps},\qquad H=2^{q-2}(\ceil K+1),\qquad m=H,\qquad \eta=\frac{2^{q-3}\eps}{s}.
\]
By Roth's theorem choose $E_1$ such that $2\sqrt{2m}\,r_3(E)\le\frac{\eps}{2s}E$ for all $E\ge E_1$, and let $T_0$ be an integer with $\eps T_0\ge E_1$ and $\frac\eps2T_0>\frac3{2\eta}+3q$. Let $T\ge T_0$ and $E\ge1$ satisfy~\eqref{eq:range}, and let $A\subseteq[1,E]$ with $\card A\ge E/q+\eps T$.

\emph{Step 1 (good core).} Since $E\ge\card A\ge\eps T\ge E_1$, \cref{lem:distinct}(b) and~\eqref{eq:range} give
\[
 \card{C_m(A)}\le\frac{\eps}{2s}E\le\frac\eps{2s}\cdot\frac{sT}{2^{q-2}}\le\frac\eps2T .
\]
Let $G=A\setminus C_m(A)$. Then $\nu_A(g)\ge m$ for all $g\in G$, and $\card G\ge E/q+\frac\eps2T$.

\emph{Step 2 (the common divisor).} Let $d=\gcd(G)$. As $G$ consists of distinct multiples of $d$ in $[1,E]$, $\card G\le E/d$; as $\card G>E/q$, we get $d<q$, hence $d\mid s$ by the definition of $q=\snd(s)$.

\emph{Step 3 (normalization).} Put $B=\{0\}\cup\{g/d:g\in G\}$ and $L=\max(G)/d$. Then $B\subseteq[0,L]$, $0,L\in B$, $\gcd(B)=1$ and $\card B=\card G+1$. Since $G\subseteq\{d,2d,\dots,dL\}$, we have $L\ge\card G\ge\frac\eps2T$; also $L\le E$. Using $L\le E$ and $T\ge2^{q-2}E/s$,
\[
 \card B>\card G\ge\frac Eq+\frac\eps2T\ge\frac Lq+\frac{\eps2^{q-2}}{2s}E\ge\Bigl(\frac1q+\eta\Bigr)L .
\]

\emph{Step 4 (the target).} Let $y=sT/d$, an integer by Step~2. Since $dL=\max(G)\le E\le sT/2^{q-2}$, we have $y\ge2^{q-2}L$; since $dL=\max(G)\ge\card G\ge\frac\eps2T$, we have $y\le sT/(\frac\eps2T)\cdot L=KL$.

\emph{Step 5.} As $L\ge\frac\eps2T_0>\frac3{2\eta}+3q$, \cref{lem:basis} gives $y\in HB$. Deleting the zero summands, $y$ is a sum of $h$ elements of $\{g/d:g\in G\}$ for some $1\le h\le H$ (and $h\ge1$ as $y>0$); multiplying by $d$, $sT=g_1+\dots+g_h$ with $g_i\in G$. Since $h\le H=m$, \cref{lem:distinct}(a) shows that $sT$ is a sum of $h$ distinct elements of $A$.
\end{proof}

\begin{proof}[Proof of \cref{cor:interval}]
Let $E=\floor{xT}$. Then $2^{q-2}E\le2^{q-2}xT\le sT$, and $\card A>(x/q+\eps)T\ge E/q+\eps T$; if $E=0$ there is no such $A$. Hence \cref{thm:inverse} applies for $T\ge T_0(s,\eps)$, which gives the first assertion and the upper bound $F(\floor{xT},sT)\le(x/q+\eps)T$. For the lower bound, if $p\nmid T$ then $q\nmid sT$ by \cref{lem:snd}(c), so the $\floor{\floor{xT}/q}\ge xT/q-2$ multiples of $q$ in $[1,\floor{xT}]$ avoid $sT$.
\end{proof}

The following proposition records the sharpness of \cref{thm:inverse}.

\begin{proposition}[Sharpness]\label{prop:sharp}
Let $s\ge1$, $q=\snd(s)$, $T\ge1$ and $t=sT$.
\begin{enumerate}[label=\textup{(\alph*)}]
\item If $q\nmid t$, the $\floor{E/q}$ multiples of $q$ in $[1,E]$ avoid $t$.
\item If $q=2$ and $E>t$, the set $[\floor{t/2}+1,E]\setminus\{t\}$ avoids $t$ and has at least $E/2+(E-t)/2-1$ elements.
\item If $q\ge3$, $q\nmid t$ and $E>t/2$, the set
\[
 A=\{x\in[1,E]:q\mid x\}\ \cup\ \{x\in(t/2,E]:\ x\not\equiv0,\ x\not\equiv t\pmod q\}
\]
avoids $t$ and has at least $E/q+\frac{q-2}q\,(E-t/2)-q$ elements.
\item If $q\ge3$, $1\le i\le q-2$, $0<\eta<1/q$ and $iE<t<(i+1)(1-1/q-\eta)E$, then the interval $[\ceil{(1-1/q-\eta)E},E]$ avoids $t$ and has at least $(1/q+\eta)E$ elements.
\end{enumerate}
Consequently, for $q\in\{2,3\}$ the range~\eqref{eq:range} in \cref{thm:inverse} cannot be replaced by $2^{q-2}E\le(1+\delta)sT$ for any $\delta>0$; and if $\theta_q$ denotes the infimum of the constants $\theta$ for which \cref{thm:inverse} holds with~\eqref{eq:range} replaced by $\theta E\le sT$, then $\theta_2=1$, $\theta_3=2$, and $\max\{2,\ q-2+1/q\}\le\theta_q\le2^{q-2}$ for $q\ge4$.
\end{proposition}

\begin{proof}
(a) is clear. (b) Two distinct elements larger than $t/2$ have sum larger than $t$, and no element equals $t$. The set has $E-\floor{t/2}-1\ge E-t/2-1$ elements.

(c) Let $\emptyset\ne B\subseteq A$ with $\sigma(B)=t$. Two elements of $A\cap(t/2,E]$ sum to more than $t$, so $B$ contains at most one element $x$ with $q\nmid x$, and every other element of $B$ is a multiple of $q$. If there is no such $x$, then $q\mid\sigma(B)=t$; otherwise $t=\sigma(B)\equiv x\not\equiv t\pmod q$. Both are impossible. For the size: there are $\floor{E/q}\ge E/q-1$ multiples of $q$; the integers in $(t/2,E]$ form an interval of length $\ell\ge E-t/2$, and since $q\nmid t$ the excluded residues $0$ and $t$ are two distinct classes, so each block of $q$ consecutive integers contributes $q-2$ elements: at least $(q-2)\floor{\ell/q}\ge\frac{q-2}q(E-t/2)-(q-2)$.

(d) Let $a=\ceil{(1-1/q-\eta)E}$. A sum of $j$ distinct elements of $[a,E]$ lies in $[ja,jE]$; if $j\le i$ it is at most $iE<t$, and if $j\ge i+1$ it is at least $(i+1)a>t$. The interval has $E-a+1\ge(1/q+\eta)E$ elements.

For the consequences, the upper bounds $\theta_q\le2^{q-2}$ are \cref{thm:inverse}. Let $\delta>0$ and put $t=sT$. If $q=2$, let $E=\floor{(1+\delta)t}$; then $E\le(1+\delta)t$, and by (b) there is a set avoiding $t$ with at least $E/2+(\delta t-1)/2-1\ge E/2+\frac{\delta s}4T$ elements for large $T$, so the conclusion of \cref{thm:inverse} fails for $\eps=\delta s/4$ and all large $T$. If $q=3$, let $E=\floor{(1+\delta)t/2}$, so $2E\le(1+\delta)t$, and take $T$ with $3\nmid T$ (then $3\nmid t$ by \cref{lem:snd}(c)); by (c) there is a set avoiding $t$ with at least $E/3+\frac13(\frac{\delta t}2-1)-3\ge E/3+\frac{\delta s}{12}T$ elements for large $T$, so the conclusion fails for $\eps=\delta s/12$ and infinitely many $T$. This proves the assertion about $q\in\{2,3\}$ and shows $\theta_2\ge1$ and $\theta_3\ge2$. The same argument with (c) shows $\theta_q\ge2$ for every $q\ge3$ (take $T$ with $p\nmid T$, where $q=p^a$). Finally let $q\ge4$ and $\theta<q-2+1/q$. Choose $\eta>0$ with $\theta':=q-2+1/q-(q-1)\eta>\max\{\theta,q-2\}$ and $\theta''\in(\max\{\theta,q-2\},\theta')$, and let $E=\ceil{sT/\theta''}$. For large $T$ we have $\theta E\le sT$ and $(q-2)E<sT<\theta'E=(q-1)(1-1/q-\eta)E$, so by (d) with $i=q-2$ there is a set of at least $E/q+\eta E\ge E/q+(\eta s/\theta'')T$ elements avoiding $sT$. Hence $\theta_q\ge q-2+1/q$.
\end{proof}

\begin{remark}\label{rem:threshold}
The exponential threshold $2^{q-2}$ comes from the dyadic iteration in \cref{lem:basis}. The example in \cref{prop:sharp}(d) shows that a bounded-basis statement for sets of density $1/q+\eta$ cannot hold below about $(q-2+1/q)L$. Structure theorems for $h$-fold sumsets, such as Lev's~\cite{Le97}, might yield a threshold linear in $q$; this would enlarge $\mathcal Q(c)$ in \cref{thm:main} accordingly. We do not pursue this here (\cref{prob:threshold}).
\end{remark}

\section{The regular range: proof of \texorpdfstring{\cref{thm:cge1}}{Theorem 1.1}}\label{sec:cge1}

\begin{proof}[Proof of \cref{thm:cge1}]
Let $1\le n\le M$ and let $A\subseteq[1,M]$ avoid $n$. Then $n\notin A$. On $A\cap[1,n-1]$ the map $x\mapsto\min\{x,n-x\}$ is injective, because two distinct elements with the same image are $x$ and $n-x$, whose sum is $n$; its image lies in $[1,\floor{n/2}]$. Therefore
\[
 \card A\le\card{A\cap[1,n-1]}+\card{A\cap[n+1,M]}\le\floor{n/2}+(M-n)=M-\ceil{n/2}.
\]
Conversely, $A^*=[\ceil{n/2},n-1]\cup[n+1,M]$ has $M-\ceil{n/2}$ elements and avoids $n$: elements larger than $n$ cannot occur in a subset with sum $n$; a single element of $[\ceil{n/2},n-1]$ is smaller than $n$; and two distinct elements of it have sum at least $2\ceil{n/2}+1>n$. Hence $F(M,n)=M-\ceil{n/2}$. If $c\ge1$ then $\floor{cn}\ge n$, so $f_c(n)=\floor{cn}-\ceil{n/2}$ and $\bigl|f_c(n)/n-(c-\frac12)\bigr|\le\frac2n$.
\end{proof}

\section{Asymptotics for \texorpdfstring{$0<c<1$}{0<c<1}}\label{sec:asymp}

\subsection{Upper bounds}

\begin{lemma}[Pairing in the middle]\label{lem:pairs}
Let $n/2<M<n$ and $D=n-M$. If $A\subseteq[1,M]$ avoids $n$, then $\card{A\cap[D,M]}\le M-n/2+1$.
\end{lemma}

\begin{proof}
The interval $[D,M]$ is invariant under $x\mapsto n-x$, and is partitioned into the pairs $\{x,n-x\}$ with $D\le x<n/2$, together with $\{n/2\}$ if $n$ is even. The set $A$ contains at most one element of each pair. There are $\ceil{n/2}-D$ pairs, so $\card{A\cap[D,M]}\le\ceil{n/2}-D+[2\mid n]\le M-n/2+1$.
\end{proof}

\begin{proposition}[Upper bounds]\label{prop:upper}
Let $\eps>0$.
\begin{enumerate}[label=\textup{(\alph*)}]
\item If $0<c\le1$, then $f_c(n)\le(c/2+\eps)n$ for all sufficiently large $n$.
\item If $0<c\le1/2$ and $q$ is a prime power with $2^{q-2}c\le1$, then $f_c(n)\le(c/q+\eps)n$ for all sufficiently large $n$ with $\snd(n)\ge q$.
\item If $1/2<c<1$ and $q$ is a prime power with $2^{q-2}(1-c)\le1$, then $f_c(n)\le\bigl(c-\frac12+\frac{1-c}q+\eps\bigr)n$ for all sufficiently large $n$ with $\snd(n)\ge q$.
\end{enumerate}
\end{proposition}

\begin{proof}
(a) Apply \cref{thm:inverse} with $s=1$ (so $q=2$), $T=n$, $E=\floor{cn}$ (so that $E\le n=sT$) and the given $\eps$: for $n\ge\max\{T_0(1,\eps),1/c\}$, every $A\subseteq[1,E]$ with $\card A\ge E/2+\eps n$ represents $n$. Hence $f_c(n)<E/2+\eps n\le(c/2+\eps)n$.

(b) Let $s=\Lambda_{q-1}$. By \cref{lem:snd}(b), $\snd(s)=q$. If $\snd(n)\ge q$ then every $k\le q-1$ divides $n$, so $s\mid n$; put $T=n/s$ and $E=\floor{cn}$. Then $2^{q-2}E\le2^{q-2}cn\le n=sT$. \Cref{thm:inverse} with $\eps s$ in place of $\eps$ shows, for $n\ge s\,T_0(s,\eps s)$ and $n\ge1/c$, that every $A\subseteq[1,E]$ with $\card A\ge E/q+\eps n$ represents $n$; hence $f_c(n)<E/q+\eps n\le(c/q+\eps)n$.

(c) Let $s,T$ be as in (b), $M=\floor{cn}$ and $D=n-M$. For $n>2/(2c-1)$ we have $n/2<M<n$. Let $A\subseteq[1,M]$ avoid $n$, and split $A=A_0\cup A_1$ with $A_0=A\cap[1,D-1]$ and $A_1=A\cap[D,M]$. By \cref{lem:pairs}, $\card{A_1}\le M-n/2+1\le(c-\frac12)n+1$. Since $D<n-cn+1$, the set $A_0$ lies in $[1,E]$ with $E=D-1\le(1-c)n$, and it avoids $n$. If $E\ge1$, then $2^{q-2}E\le2^{q-2}(1-c)n\le sT$, and \cref{thm:inverse} with $\eps s/2$ in place of $\eps$ gives, for $n\ge s\,T_0(s,\eps s/2)$, that $\card{A_0}<E/q+\frac\eps2n\le\bigl(\frac{1-c}q+\frac\eps2\bigr)n$. Hence $\card A\le\bigl(c-\frac12+\frac{1-c}q+\frac\eps2\bigr)n+1$, which is at most $\bigl(c-\frac12+\frac{1-c}q+\eps\bigr)n$ for $n\ge2/\eps$.
\end{proof}

\subsection{Constructions}

\begin{proposition}[Constructions]\label{prop:constructions}
Let $n\ge1$ and $M\ge0$ be integers.
\begin{enumerate}[label=\textup{(\alph*)}]
\item If $q\ge2$ and $q\nmid n$, then $Q=\{x\in[1,M]:q\mid x\}$ avoids $n$, and $\card Q=\floor{M/q}$.
\item If $q\ge3$, $q\nmid n$ and $M<n$, then $Q\cup U$ avoids $n$, where
\[
 U=\bigl\{x\in[\ceil{n/2},M]:\ x\not\equiv0,\ x\not\equiv n\pmod q\bigr\},
\]
and $\card{Q\cup U}\ge M/q+\frac{q-2}q\,(M-n/2)-(q-1)$.
\item For every $k\ge2$ the set $\{x\in\N:\ n/(k+1)<x\le n/k\}$ avoids $n$, and so do $[\ceil{n/2},n-1]$ and $[1,k]$ whenever $k(k+1)/2<n$.
\item If $M<n$ and $D:=n-M$ satisfies $n/4<D\le(n+2)/3$, then $\{x\in\N:\ n/4<x<D\}\cup[\ceil{n/2},M]$ avoids $n$ and has at least $n/4-1$ elements.
\end{enumerate}
\end{proposition}

\begin{proof}
(a) Every subset sum of $Q$ is divisible by $q$.

(b) Let $\emptyset\ne B\subseteq Q\cup U$ with $\sigma(B)=n$. Two distinct elements of $U$ have sum at least $2\ceil{n/2}+1>n$, so $\card{B\cap U}\le1$. If $B\cap U=\emptyset$ then $q\mid\sigma(B)=n$; if $B\cap U=\{x\}$ then $n=\sigma(B)\equiv x\not\equiv n\pmod q$. Both are impossible. The sets $Q$ and $U$ are disjoint, and $\card Q=\floor{M/q}\ge M/q-\frac{q-1}q$. If $M<\ceil{n/2}$ then $M-n/2<0$ and the bound follows from $\card Q\ge M/q-\frac{q-1}{q}$. If $M\ge\ceil{n/2}$, the interval $[\ceil{n/2},M]$ has $\ell=M-\ceil{n/2}+1\ge M-n/2+\frac12$ elements and contains $\floor{\ell/q}\ge(\ell-q+1)/q$ disjoint blocks of $q$ consecutive integers; since $q\nmid n$, the residues $0$ and $n$ are distinct, so each block contains $q-2$ elements of $U$. Hence $\card U\ge\frac{q-2}q(M-n/2)-\frac{(q-2)(q-3/2)}q$, and the claim follows from $\frac{q-1}q+\frac{(q-2)(q-3/2)}q=q-\frac52+\frac2q\le q-1$.

(c) For the first set: a sum of at most $k$ distinct elements, each at most $n/k$, is smaller than $n$ (as $k\ge2$), and a sum of at least $k+1$ elements, each larger than $n/(k+1)$, exceeds $n$. The interval $[\ceil{n/2},n-1]$ was treated in the proof of \cref{thm:cge1}. All subset sums of $[1,k]$ are at most $k(k+1)/2<n$.

(d) Let $B$ be a subset with sum $n$. It contains at most one element of $[\ceil{n/2},M]$, since two of them sum to more than $n$. If it contains none, it consists of elements of $(n/4,D)$, all at most $D-1$: three distinct ones sum to at most $3D-6\le n-4<n$, and four of them sum to more than $n$. If it contains one element $y\in[\ceil{n/2},M]$, the remaining elements have sum $n-y\in[D,\floor{n/2}]$; but a single element of $(n/4,D)$ is smaller than $D$, and two of them sum to more than $n/2$. The size is at least $(D-1-n/4)+(M-n/2+\frac12)=n/4-\frac12$.
\end{proof}

\subsection{Proofs of \texorpdfstring{\cref{thm:main}}{Theorem 1.4} and its corollaries}

\begin{proof}[Proof of \cref{thm:main}]
(i) The upper bound is \cref{prop:upper}(a). For odd $n$, \cref{prop:constructions}(a) with $q=2$ and $M=\floor{cn}$ gives $f_c(n)\ge\floor{\floor{cn}/2}\ge(cn-1)/2-\frac12=cn/2-1$. Together these give $c/2-1/n\le f_c(n)/n\le c/2+\eps$ for all large odd $n$, so the limit along odd $n$ is $c/2$; since the upper bound holds for all $n$, the $\limsup$ over all $n$ is also $c/2$.

(ii) Let $q\in\mathcal Q(c)$, so $2^{q-2}\min\{c,1-c\}\le1$. If $c\le1/2$, then $\varphi_q(c)=c/q$ and $2^{q-2}c\le1$, so \cref{prop:upper}(b) gives the upper bound along $\snd(n)\ge q$; and if $\snd(n)=q$ then $q\nmid n$, and \cref{prop:constructions}(a) gives $f_c(n)\ge\floor{\floor{cn}/q}\ge cn/q-2\ge\varphi_q(c)n-q$. If $c>1/2$, then $\varphi_q(c)=c-\frac12+\frac{1-c}q$ and $2^{q-2}(1-c)\le1$, so \cref{prop:upper}(c) gives the upper bound; and if $\snd(n)=q$, then $q\nmid n$, $M=\floor{cn}<n$, and \cref{prop:constructions}(b) (for $q\ge3$; for $q=2$ use (a)) gives
\begin{align*}
 f_c(n)&\ge\frac Mq+\frac{q-2}q\Bigl(M-\frac n2\Bigr)-(q-1)
 \ \ge\ \frac{cn-1}q+\frac{q-2}q\Bigl(cn-1-\frac n2\Bigr)-(q-1)\\
 &=\varphi_q(c)\,n-\frac{q-1}q-(q-1),
\end{align*}
because $\frac cq+\frac{q-2}q(c-\frac12)=c-\frac12+\frac{1-c}q$; the right-hand side is at least $\varphi_q(c)\,n-q$. Since the class $\{n:\snd(n)=q\}$ is infinite (it contains $\Lambda_{q-1}(pk+1)$ for all $k\ge0$, where $q=p^a$), the two bounds give both displayed limits.
\end{proof}

\begin{proof}[Proof of \cref{cor:irregular}]
For fixed $c\in(0,1)$, both $q\mapsto c/q$ and $q\mapsto c-\frac12+\frac{1-c}q$ are strictly decreasing, hence so is $q\mapsto\varphi_q(c)$; therefore the limits $\varphi_q(c)$, $q\in\mathcal Q(c)$, are pairwise distinct limit points of $f_c(n)/n$ by \cref{thm:main}(ii). Since $2,3\in\mathcal Q(c)$, the sequence has at least two limit points and does not converge. A prime power $q$ belongs to $\mathcal Q(c)$ if and only if $q\le2+\log_2(1/\min\{c,1-c\})$, which gives the lower bound for the number of limit points. The last assertion follows from \cref{thm:cge1}.
\end{proof}

\begin{proof}[Proof of \cref{cor:half}]
The odd case is \cref{thm:main}(i). Let $n$ be even and $E=n/2=\floor{n/2}$. \Cref{thm:inverse} with $s=2$ (so $q=3$), $T=n/2$ and $E=T$ (so $2E=sT$) shows that $f_{1/2}(n)<E/3+\eps T=(1/6+\eps/2)n$ for large even $n$. Conversely, if $3\nmid n$ the multiples of $3$ give $f_{1/2}(n)\ge\floor{n/6}$, and if $3\mid n$ the set $\{x:n/3<x\le n/2\}$ (\cref{prop:constructions}(c) with $k=2$) has $n/6$ elements and avoids $n$. Hence $f_{1/2}(n)/n\to1/6$ along even $n$. Every limit point is a limit along a subsequence lying eventually in the odd or in the even integers, so the set of limit points is $\{1/6,1/4\}$.
\end{proof}

The even case of \cref{cor:half} is, in substance, Alon's $f(N,2N)=(1/3+o(1))N$~\cite[Corollary~2.6]{Al87}.

\subsection{Remarks}

\begin{remark}[Lower bounds for all $n$]\label{rem:liminf}
Let $0<c\le1$, $n\ge1$ and $M=\floor{cn}$. Intersecting the sets of \cref{prop:constructions}(c) with $[1,M]$, and using $[\ceil{n/2},\min\{M,n-1\}]$ when $k=1$, we obtain for every $k\ge1$
\[
 f_c(n)\ \ge\ \Bigl(\min\Bigl\{c,\frac1k\Bigr\}-\frac1{k+1}\Bigr)n-1 .
\]
Hence, for $0<c<1$ and $t=\ceil{1/c}$, $\liminf_{n}f_c(n)/n\ge\max\{c-\frac1t,\ \frac1{t(t+1)}\}>0$; this was observed in the discussion on~\cite{Bl}. For $c=1/2$ it gives $1/6$, which is the true value of the $\liminf$ by \cref{cor:half}. For $2/3\le c<3/4$, \cref{prop:constructions}(d) improves it to $\liminf_n f_c(n)/n\ge1/4$.
\end{remark}

\begin{remark}[On the hypothesis $q\in\mathcal Q(c)$]\label{rem:intervals}
Let $2/3\le c<3/4$. Then $D=n-\floor{cn}$ satisfies $D\ge(1-c)n>n/4$ and $D<(1-c)n+1\le n/3+1$, hence $D\le\ceil{n/3}\le(n+2)/3$; so \cref{prop:constructions}(d) gives $f_c(n)\ge n/4-1$ for every $n$. Since $\varphi_q(c)=c-\frac12+\frac{1-c}q<\frac14$ as soon as $q>(1-c)/(\frac34-c)$, the conclusion of \cref{thm:main}(ii) fails for all such $q$. For instance, for $c=2/3$ it fails for every $q\ge5$ (here $\mathcal Q(2/3)=\{2,3\}$, and $\varphi_4(2/3)=1/4$). Similarly, for $c\le1/2$ the interval $(n/(k+1),cn]$ with $k=\floor{1/c}$ beats the multiples of $q$ as soon as $q$ is large compared with $1/c$.
\end{remark}

\begin{remark}[Small $n$]
\Cref{prop:constructions}(c) also gives $f_c(n)\ge\min\{\floor{cn},\floor{\sqrt{2n}-1/2}\}$, which dominates for small $n$ when $c$ is small; this explains the numerical data of \cref{sec:numerics} for $c\le1/2$.
\end{remark}

\section{A zero-sum proof of non-convergence}\label{sec:zerosum}

This section gives a second, independent proof that $f_c(n)/n$ diverges for $0<c<1$, based on bounded zero-sums. It yields weaker constants than \cref{thm:main}; it is the argument followed by the Lean proof of the statement \lean{erdos361_irregular} (\cref{sec:lean}).

\begin{theorem}[Alon~{\cite[Theorem~1.1]{Al87}}]\label{thm:alon}
For every integer $k\ge2$ and every $\eps>0$ there is $N_0$ such that, for every $N\ge N_0$, every $X\subseteq\Z/N\Z$ with $\card X>(1/k+\eps)N$ contains a nonempty $Y$ with $\card Y\le k$ and $\sum_{y\in Y}y=0$.
\end{theorem}

For prime moduli the theorem follows at once from the Dias da Silva--Hamidoune theorem~\cite{DSH94} (see also~\cite{ANR96}): if $p$ is prime and $X\subseteq\Z/p\Z$, then the set $k^\wedge X$ of sums of $k$ distinct elements of $X$ satisfies $\card{k^\wedge X}\ge\min\{p,\ k\card X-k^2+1\}$.

\begin{lemma}\label{lem:dsh}
Let $p$ be prime, $k\ge1$, and $X\subseteq\Z/p\Z$ with $\card X\ge p/k+k$. Then there is $Y\subseteq X$ with $\card Y=k$ and $\sum_{y\in Y}y=0$. In particular \cref{thm:alon} holds for prime $N\ge k^2/\eps$.
\end{lemma}

\begin{proof}
$k\card X-k^2+1\ge p+1$, so $k^\wedge X=\Z/p\Z\ni0$. If $\card X>(1/k+\eps)p$ and $p\ge k^2/\eps$, then $\card X\ge p/k+k$.
\end{proof}

\begin{proposition}[Two zero-sums; after {\cite[Corollary~2.6]{Al87}}]\label{prop:twozero}
Fix $\rho\ge2$ and $\delta>0$. For all sufficiently large $E$ the following holds uniformly for all even $\tau$ with $2E\le\tau\le\rho E$: if $A\subseteq[1,E]$ and $\tau\notin\Sig(A)$, then
\[
 \card A\le\frac{\tau}{2K}+\delta E,\qquad K=\Bigl\lfloor\frac{3\tau}{2E}\Bigr\rfloor .
\]
\end{proposition}

\begin{proof}
Put $\mu=\tau/2$, so $E\le\mu\le\rho E/2$ and $3\le K\le3\rho/2$; thus $K$ ranges over a fixed finite set. Apply \cref{thm:alon} with $\eps=\delta/\rho$ for each such $K$, and take $E$ so large that $E\ge N_0$ for all these $K$ and $K<\delta E/2$.

Suppose $\card A>\mu/K+\delta E$. Reduction modulo $\mu$ is injective on $[1,E]$, and $\delta E\ge2\delta\mu/\rho$, so the image of $A$ has more than $(1/K+\delta/\rho)\mu$ elements. Hence there is a nonempty $B_1\subseteq A$ with $\card{B_1}\le K$ and $\mu\mid\sigma(B_1)$. Moreover $0<\sigma(B_1)\le KE\le3\mu$, and $\sigma(B_1)=3\mu$ would force $\card{B_1}=K\ge3$ with all elements equal to $E$, which is impossible. So $\sigma(B_1)\in\{\mu,2\mu\}$. If $\sigma(B_1)=2\mu=\tau$ we have a contradiction. Otherwise $\sigma(B_1)=\mu$, and
\[
 \card{A\setminus B_1}>\frac\mu K+\delta E-K>\frac\mu K+\frac{\delta E}2\ge\Bigl(\frac1K+\frac\delta\rho\Bigr)\mu ,
\]
so a second application gives a nonempty $B_2\subseteq A\setminus B_1$ with $\sigma(B_2)\in\{\mu,2\mu\}$. Then $B_2$ or $B_1\cup B_2$ has sum $\tau$, a contradiction.
\end{proof}

\begin{theorem}\label{thm:zerosum-irregular}
For every $0<c<1$, $\limsup_{n\ \mathrm{even}}f_c(n)/n<c/2\le\liminf_{n\ \mathrm{odd}}f_c(n)/n$. In particular $f_c(n)/n$ does not converge.
\end{theorem}

\begin{proof}
The inequality for odd $n$ is \cref{thm:main}(i) (the lower bound there is elementary). Suppose first $c\le1/2$ and let $J=\floor{3/(2c)}$. Since $J>3/(2c)-1\ge1/c$, we have $1/(2J)<c/2$; choose $\delta>0$ with $\beta=1/(2J)+\delta c<c/2$. For large even $n$ put $M=\floor{cn}$; then $2M\le n\le(2/c)M$, and \cref{prop:twozero} with $E=M$, $\tau=n$, $\rho=2/c$ gives $f_c(n)\le n/(2K)+\delta M$ with $K=\floor{3n/(2M)}\ge J$. Hence $f_c(n)/n\le\beta$ for all large even $n$.

Now let $1/2<c<1$, $d=1-c$ and $J=\floor{3/(2d)}$; as $d<1/2$, again $J>1/d$. Let $g=d/2-1/(2J)>0$, $\delta=g/(3d)$, $\eta=g/3$, so that $c-\frac12+\frac1{2J}+\delta d+\eta=\frac c2-\frac g3$. For a large even $n$ let $M=\floor{cn}$, $D=n-M$, and let $A\subseteq[1,M]$ avoid $n$. By \cref{lem:pairs}, $\card{A\cap[D,M]}\le(c-\frac12)n+1$. Since $dn\le D<dn+1$, for large $n$ we have $2(D-1)\le n\le(2/d)(D-1)$, and \cref{prop:twozero} with $E=D-1$, $\tau=n$, $\rho=2/d$ gives $\card{A\cap[1,D-1]}\le n/(2K)+\delta dn$ with $K=\floor{3n/(2(D-1))}\ge J$. Hence $f_c(n)/n\le c-\frac12+\frac1{2J}+\delta d+\frac1n\le\frac c2-\frac g3$ for $n\ge1/\eta$.
\end{proof}

\begin{remark}
If $\tau=2p$ with $p$ prime, the modulus $\mu$ in the proof of \cref{prop:twozero} is prime, and \cref{lem:dsh} can replace \cref{thm:alon}. Restricting the even subsequence in the proof of \cref{thm:zerosum-irregular} to $n=2p$ is harmless, since there are infinitely many primes. This is how the formal proof avoids the general case of Alon's theorem: the Dias da Silva--Hamidoune bound is derived there from the Combinatorial Nullstellensatz~\cite{Al99} through the Alon--Nathanson--Ruzsa polynomial method~\cite{ANR96}.
\end{remark}

\section{Numerical evidence}\label{sec:numerics}

We computed $F(\floor{cn},n)$ exactly for the twelve values $c\in\{\frac15,\frac3{10},\frac13,\frac25,\frac9{20},\frac12,\frac35,\frac23,\frac7{10},\frac34,\frac45,\frac9{10}\}$ and all $40\le n\le127$ ($1056$ instances) by a depth-first branch-and-bound search. The elements $x=1,2,\dots,M$ are decided in increasing order; the set $R\subseteq[0,n-1]$ of subset sums of the chosen elements (including $0$) is stored as a $128$-bit mask, and $x$ may be added if and only if $n-x\notin R$. At a node at which $1,\dots,i$ have been decided, the search is pruned using the bound
\[
 \card{A_{\mathrm{cur}}}+\#\{x\in(i,M]:n-x\notin R\}-\Bigl\lfloor\tfrac12\#V\Bigr\rfloor,
\]
where $V$ is the set of $x\in[\max(i+1,n-M),\min(M,n-i-1)]$, $x\ne n/2$, with $x\notin R$ and $n-x\notin R$. The set $V$ is invariant under $x\mapsto n-x$, both elements of each pair $\{x,n-x\}\subseteq V$ are individually addable, and at most one of them can be added; this justifies the bound. The results were cross-checked against an independent implementation (descending order, without the pair bound) on all $300$ instances with $n\le64$, with no discrepancy. The constructions were checked by a subset-sum bitset computation: \cref{prop:constructions}(b) with $q=\snd(n)$ for all $n<300$ and $1\le M<n$ ($44{,}551$ instances, avoidance), and its size bound for $q\in\{3,4,5,7,8,9,11\}$, $q\nmid n$, $n<800$, $1\le M<n$ ($1{,}833{,}234$ instances); \cref{prop:constructions}(c) for $n<400$ and $2\le k\le11$; \cref{prop:constructions}(d) for all admissible $n<700$ ($20{,}880$ instances); and \cref{prop:sharp}(b)--(d) for $t<160$ ($185{,}741$ instances).

\emph{Forum values.} For $c=3/4$ we obtain $F(\floor{3n/4},n)=34,37,32,38,35$ for $n=100,\dots,104$, in agreement with the values reported in the discussion of the problem on~\cite{Bl}. The extremal sets found are, respectively: $Q\cup U$ with $q=3$; the even numbers; $[1,6]\cup[51,76]$; the even numbers; and $Q\cup U$ with $q=3$.

\emph{Exactness for odd $n$.} Let $k_n$ be the largest integer $k$ with $k(k+1)/2<n$. For all $528$ computed instances with $n$ odd we found
\[
 f_c(n)=\max\bigl\{\floor{\floor{cn}/2},\ \min\{\floor{cn},k_n\}\bigr\},
\]
the two terms coming from the even numbers and from $[1,k_n]$. In particular $f_c(n)=\floor{\floor{cn}/2}$ for every $c\ge2/5$ in our list and every odd $n\in[41,127]$, and for $c=1/3$ (respectively $c=3/10$) for every odd $n\ge61$ (respectively $n\ge81$).

\begin{figure}[ht]
\centering
\includegraphics[width=\linewidth]{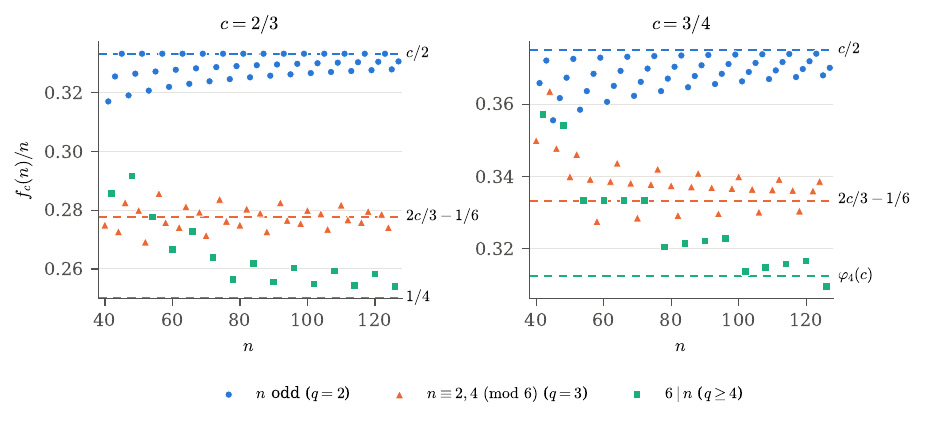}
\caption{Exact values of $f_c(n)/n$ for $40\le n\le127$, for $c=2/3$ (left) and $c=3/4$ (right), by class of $n$. Dashed lines: the limits $c/2$ (odd $n$) and $2c/3-1/6$ ($n\equiv2,4\bmod6$) of \cref{thm:main}; for $c=3/4$ also $\varphi_4(3/4)=5/16$, the limit along $n\equiv6\pmod{12}$; for $c=2/3$ the lower bound $1/4$ of \cref{rem:intervals}, valid for all $n$.}\label{fig:fc}
\end{figure}

\emph{Comparison with \cref{thm:main}.} \Cref{tab:numerics} lists, for each $c$, the mean of $f_c(n)/n$ over the $n\in[90,127]$ with $\snd(n)=2$ ($19$ values), $\snd(n)=3$ ($12$ values) and $\snd(n)=4$ ($n=90,102,114,126$), next to $\varphi_q(c)$ when $q\in\mathcal Q(c)$. For $c\ge3/5$ every listed mean differs from the corresponding limit by less than $0.016$, consistent with lower-order terms, and \cref{fig:fc} shows the convergence along the three classes for $c=2/3$ and $c=3/4$. For $c\le1/2$ the range $n\le127$ is too small: the sets $[1,k]$ with $k\approx\sqrt{2n}$ have about $0.12\,n$ elements in this range, which exceeds the asymptotic values.

\begin{table}[ht]
\centering
\caption{Mean of $f_c(n)/n$ over $90\le n\le127$ in the classes $\snd(n)=q$, and the limit $\varphi_q(c)$ of \cref{thm:main} when $q\in\mathcal Q(c)$ (otherwise ``---'').}\label{tab:numerics}
\begin{tabular}{@{}lcccccc@{}}
\toprule
& \multicolumn{2}{c}{$q=2$ ($n$ odd)} & \multicolumn{2}{c}{$q=3$} & \multicolumn{2}{c}{$q=4$}\\
\cmidrule(lr){2-3}\cmidrule(lr){4-5}\cmidrule(l){6-7}
$c$ & mean & $\varphi_2(c)$ & mean & $\varphi_3(c)$ & mean & $\varphi_4(c)$\\
\midrule
$1/5$ & 0.1266 & 0.1000 & 0.1279 & 0.0667 & 0.1257 & 0.0500 \\
$3/10$ & 0.1452 & 0.1500 & 0.1279 & 0.1000 & 0.1257 & --- \\
$1/3$ & 0.1622 & 0.1667 & 0.1279 & 0.1111 & 0.1257 & --- \\
$2/5$ & 0.1965 & 0.2000 & 0.1382 & 0.1333 & 0.1257 & --- \\
$9/20$ & 0.2201 & 0.2250 & 0.1641 & 0.1500 & 0.1581 & --- \\
$1/2$ & 0.2452 & 0.2500 & 0.1807 & 0.1667 & 0.1761 & --- \\
$3/5$ & 0.2954 & 0.3000 & 0.2347 & 0.2333 & 0.2129 & --- \\
$2/3$ & 0.3303 & 0.3333 & 0.2778 & 0.2778 & 0.2547 & --- \\
$7/10$ & 0.3452 & 0.3500 & 0.3001 & 0.3000 & 0.2755 & --- \\
$3/4$ & 0.3703 & 0.3750 & 0.3356 & 0.3333 & 0.3153 & 0.3125 \\
$4/5$ & 0.3964 & 0.4000 & 0.3681 & 0.3667 & 0.3571 & 0.3500 \\
$9/10$ & 0.4450 & 0.4500 & 0.4391 & 0.4333 & 0.4402 & 0.4250 \\
\bottomrule
\end{tabular}
\end{table}

\section{Formal verification}\label{sec:lean}

The formalization, available at \url{https://github.com/matejasevick/erdos361}, uses Lean~4 (toolchain \lean{v4.31.0}) and Mathlib (revision \lean{fabf563a}). The trusted statement layer is a file that imports only Mathlib and consists of four definitions:
\begin{align*}
 \lean{Avoids A n} &\iff \forall B\subseteq A,\ \textstyle\sum_{b\in B}b=n\ \Rightarrow\ B=\emptyset,\\
 \lean{Avoiders M n}&=\{A\subseteq[1,M]:\ \lean{Avoids A n}\},\qquad
 \lean{F M n}=\max_{A\in\lean{Avoiders M n}}\card A,\\
 \lean{Fc c n}&=\lean{F}\,\floor{c\,n}_+\,n .
\end{align*}
For $n\ge1$, $\lean{Avoids A n}$ is equivalent to $n\notin\Sig(A)$, $\lean{F}$ is the function $F$ of the introduction (the maximum exists since $\emptyset$ avoids $n$), and $\lean{Fc c n}=f_c(n)$ because $\floor{\cdot}_+$, the natural-number floor, agrees with $\floor\cdot$ on $[0,\infty)$. A file \lean{Challenge.lean}, which imports only Mathlib, repeats these four definitions verbatim and states the theorems without proofs; \lean{Solution.lean} proves each of them by a direct term assignment from the development. The repository also contains configurations for the tool Comparator~\cite{Comparator}, which checks mechanically that proved statements coincide with challenge statements and that only permitted axioms are used.

\begin{table}[ht]
\centering
\caption{Lean statements (namespace \lean{Erdos361.Statement}) and the corresponding results. The statements above the rule are proved in the earlier development, those below it in the module \lean{Asymptotics}. All eleven have been machine-checked.}\label{tab:lean}
\small
\begin{tabular}{@{}lll@{}}
\toprule
Lean name & Result & Module\\
\midrule
\lean{erdos361_cge1} & \cref{thm:cge1} ($F(M,n)$ formula) & \lean{Core}\\
\lean{erdos361_irregular} & non-convergence, $0<c<1$ & \lean{Core}\\
\lean{inverse_s2_restricted} & \cref{thm:inverse}, $s=2$, $2E<sT<3E$ & \lean{InverseMain}\\
\midrule
\lean{inverse_s1}, \lean{inverse_s2} & \cref{thm:inverse}, $q\in\{2,3\}$ & \lean{Asymptotics}\\
\lean{erdos361_upper_half} & \cref{thm:main}(i), upper bound & \lean{Asymptotics}\\
\lean{erdos361_odd_limit} & \cref{thm:main}(i), odd $n$ & \lean{Asymptotics}\\
\lean{erdos361_mod6_limit_low} & \cref{thm:main}(ii), $q=3$, $c\le\frac12$ & \lean{Asymptotics}\\
\lean{erdos361_mod6_limit_high} & \cref{thm:main}(ii), $q=3$, $c>\frac12$ & \lean{Asymptotics}\\
\lean{erdos361_cge1_limit} & \cref{thm:cge1}, limit form & \lean{Asymptotics}\\
\lean{erdos361_two_limit_points} & \cref{cor:irregular}, $q\in\{2,3\}$ & \lean{Asymptotics}\\
\bottomrule
\end{tabular}
\end{table}

\emph{The earlier development.} The module \lean{Core} proves the formula of \cref{thm:cge1} and the non-convergence; following \cref{sec:zerosum}, it derives the Dias da Silva--Hamidoune bound from Mathlib's Combinatorial Nullstellensatz and works along $n=2p$. The modules \lean{GoodCore}, \lean{Covering}, \lean{Freiman}, \lean{InverseGlue} and \lean{InverseMain} prove the case $s=2$, $2E<sT<3E$, $3\nmid sT$ of \cref{thm:inverse}, in the form $\card A\ge(1/3+\eps)E$ (the hypothesis $3\nmid sT$ is superfluous), including a proof of \cref{thm:freiman} from scratch (following the large-diameter induction in~\cite[Chapter~1]{Na96}), a version of \cref{lem:distinct} and a covering lemma.

\emph{The new module.} The module \lean{Asymptotics} formalizes \cref{lem:distinct} with an arbitrary threshold $m$, the zero-padded covering argument of \cref{lem:basis} for $q\in\{2,3\}$, \cref{thm:inverse} for $q\in\{2,3\}$ in the form ``$\card A\ge(1/q+\eps)E$ implies $t\in\Sig(A)$ for all $t\in[2^{q-2}E,KE]$ (with $2\mid t$ when $q=3$)'', the constructions of \cref{prop:constructions}(a),(b) for $q\in\{2,3\}$, and the limits of \cref{thm:main} along $\snd(n)=q$ for $q\in\{2,3\}$. For $s$ with $\snd(s)=q\in\{2,3\}$, this form with $t=sT$ yields the representation statement of \cref{thm:inverse}: since $E\ge\card A\ge\eps T$ and $sT\ge E$, we have $sT\le(s/\eps)E$ and $\card A\ge(1/q+\eps/s)E$. This short deduction and the bound on the number of summands are not part of the formalization.

\emph{Verification.} The development consists of nine modules with about $5000$ lines of Lean code, of which about $1000$ are in \lean{Asymptotics}. On 11--12~November~2026 the repository was compiled from source with Lean~4.31.0 against Mathlib at commit \lean{fabf563a} (tag \lean{v4.31.0}), one module at a time; its file \lean{VERIFICATION.md} records the exit code of every step, the axiom reports and the SHA-256 hashes of the compiled files. All modules, \lean{Challenge.lean} and \lean{Solution.lean} compiled without errors, so the Lean kernel accepted every declaration, and for each of the eleven statements of \cref{tab:lean} the command \texttt{\#print axioms} reports exactly \lean{[propext, Classical.choice, Quot.sound]}. The development declares no axioms, instances, coercions or notation, and does not use \lean{sorry}, \lean{native_decide}, \lean{unsafe} or \lean{extern}; the only occurrences of \lean{sorry} in the repository are the eleven statements of \lean{Challenge.lean}. After removal of comments, each statement in \lean{Solution.lean} is character-for-character identical to the corresponding statement in \lean{Challenge.lean}, and the definitions in \lean{Challenge.lean} are identical to those in the statement layer. The Comparator configurations were not run as part of this verification.

\section{Open problems}\label{sec:open}

\begin{problem}\label{prob:classes}
Let $0<c<1$ and let $q\notin\mathcal Q(c)$ be a prime power. Does $f_c(n)/n$ converge along $\snd(n)=q$, and if so, to what? The smallest open cases are $q=4$ (that is, $n\equiv6\pmod{12}$) for $1/4<c<3/4$; for $c=2/3$, \cref{rem:intervals} and \cref{thm:main} give $1/4\le\liminf\le\limsup\le\varphi_3(2/3)=5/18$ along this class, and the data in \cref{tab:numerics} suggest the value $1/4$.
\end{problem}

\begin{problem}\label{prob:threshold}
Determine the optimal range constant $\theta_q$ in \cref{thm:inverse} (see \cref{prop:sharp}). Is $\theta_q=O(q)$? Is $\theta_q=\max\{2,q-2+1/q\}$?
\end{problem}

\begin{problem}
Determine $\liminf_{n\to\infty}f_c(n)/n$ and the set of limit points of $f_c(n)/n$ for $0<c<1$, $c\ne1/2$. Is the set of limit points finite?
\end{problem}

\begin{conjecture}
For every $0<c\le1$ there is $n_0(c)$ such that $f_c(n)=\floor{\floor{cn}/2}$ for all odd $n\ge n_0(c)$.
\end{conjecture}

The conjecture is supported by the computations of \cref{sec:numerics}. Proving it would require an exact, rather than asymptotic, version of the case $s=1$ of \cref{thm:inverse}. In the same spirit, the discussion on~\cite{Bl} contains the assertion $f_c(n)\le cn/2+2$ for $1/2<c<1$; all $528$ computed values with $c>1/2$ satisfy $f_c(n)\le\floor{cn}/2+2$.

\section*{Acknowledgments}
This paper is based on ideas produced with a research pipeline combining the AI systems GPT-5 (OpenAI) and Claude (Anthropic). The literature search, the numerical computations, and the Lean formalization were prepared with the assistance of Claude. The author thanks the participants of the online discussion of Problem~\#361 for the constructions and data recalled in \cref{subsec:literature,sec:numerics}. The author takes full responsibility for the content of this paper. 

\section*{Data availability}
The repository \url{https://github.com/matejasevick/erdos361} contains the Lean sources and the verification record described in \cref{sec:lean}. The C program used for \cref{sec:numerics}, the computed values and the build logs are available from the author.

\end{document}